\documentclass[12pt]{article}

\usepackage[centertags]{amsmath}
\usepackage{amssymb}
\usepackage{amsthm}
\usepackage{verbatim}
\usepackage[dvips]{graphics}
\usepackage{enumerate}

\newtheorem{proposition}{Proposition}
\newtheorem{lemma}[proposition]{Lemma}
\newtheorem{definition}[proposition]{Definition}

\newtheorem{theorem}[proposition]{Theorem}

\newcommand{\remark}{\noindent\textbf{Remark}}

\newcommand{\dif}{\mathrm{d}}
\newcommand{\eps}{\varepsilon}

\newcommand{\tldvphi}{\tilde{\varphi}}
\newcommand{\tldh}{\tilde{h}}

\newcommand{\vphi}{\varphi}

\newcommand{\CC}{\mathbb{C}}
\newcommand{\DD}{\mathbb{D}}
\newcommand{\RR}{\mathbb{R}}
\newcommand{\Fl}{\mathcal{F}}

\newcommand{\dist}{\mathrm{dist}}
\newcommand{\inrad}{\mathrm{inrad}}

\newcommand{\Area}{\mathrm{Area}}
\newcommand{\Vol}{\mathrm{Vol}}
\newcommand{\vol}{\mathrm{vol}}

\newcommand{\supp}{\mathrm{supp}}

\newcommand{\CAP}{\mathrm{cap}}
\newcommand{\hole}{\mathrm{hole}}
\newcommand{\pr}{\mathrm{pr}}

\begin{document}

\title{On the Inner Radius of Nodal Domains}
\author{Dan Mangoubi}
\date{}
\maketitle
\abstract
 Let $M$ be a closed Riemannian manifold.
 We consider the inner radius of a nodal domain for a large eigenvalue $\lambda$.
 We give upper and lower bounds on the inner radius of the type
 $C/\lambda^\alpha(\log\lambda)^\beta$. Our proof is based on
 a local behavior of eigenfunctions discovered by Donnelly and
 Fefferman and a Poincar\'{e} type inequality proved by Maz'ya.
 Sharp lower bounds are known
 only in dimension two. We give an account of this case too.
\vspace{1.5ex}

 \vspace{1ex}
\section{Introduction and Main Results}
Let $(M, g)$ be a closed Riemannian manifold of dimension $n$. Let
$\Delta$ be the Laplace--Beltrami operator on $M$. Let
$0<\lambda_1\leq \lambda_2 \leq \ldots$ be the eigenvalues of
$\Delta$.
 Let $\vphi_\lambda$ be an eigenfunction of $\Delta$ with eigenvalue $\lambda$.
A \emph{nodal domain} is a connected component of
$\{\vphi_\lambda\neq 0\}$.

We are interested in the asymptotic geometry of the nodal domains.
In particular, in this paper we consider the inner radius of nodal
domains.

 Let $r_\lambda$ be the
inner radius of the $\lambda$-nodal domain $U_\lambda$. Let $C_1,
C_2, \ldots$ denote constants which depend only on $(M, g)$.
 We prove
\begin{theorem}
\label{thm:inrad}
  Let $M$ be a closed Riemannian manifold of dimension $n\geq 3$.
  Then
  $$ \frac{C_1}{\sqrt{\lambda}} \geq r_\lambda \geq
  \frac{C_2}{\lambda^{k(n)}(\log{\lambda})^{2n-4}}\, ,$$
  where $k(n)=n^2-15n/8+1/4$.
\end{theorem}

In dimension two we have the following sharp bound
\begin{theorem}
\label{thm:inrad2} Let $\Sigma$ be a closed Riemannian surface.
Then
\begin{equation*}
  \frac{C_3}{\sqrt{\lambda}}\geq r_\lambda \geq
  \frac{C_4}{\sqrt{\lambda}} \, .
\end{equation*}
\end{theorem}
%
\subsection{Upper Bound}
We remark that the upper bound is more or less standard and has been used
in the literature~(e.g.~\cite{don-fef-yau}). However, we explain it here also.

 We observe that $\lambda=\lambda_1(U_\lambda)$. This is true since the
$\lambda$-eigenfunction does not vanish in $U_\lambda$
(\cite{chavel}, ch.~I.5). Therefore, the existence of the upper
bound in Theorems~\ref{thm:inrad} and~\ref{thm:inrad2} follows
from the following general upper bound on $\lambda_1$ of domains
$\Omega\subseteq M$.
\begin{theorem}
  \label{thm:inrad-gen}
  $$ \lambda_1(\Omega) \leq \frac{C_5}{\inrad(\Omega)^2}\, .$$
\end{theorem}
The proof of this theorem is given in \S\ref{sec:gen-bounds}.
\subsection{Lower Bound}
For the lower bound on the inner radius in dimensions $\geq 3$, we
give a proof in \S\ref{sec:lower3} which is based on a local
behavior of eigenfunctions discovered by H.~Donnelly and
C.~Fefferman (Theorem~\ref{thm:loc-cour}). The same proof gives in
dimension two the bound $C/\sqrt{\lambda\log\lambda}$.

In order to get rid of the factor $\sqrt{\log\lambda}$ in
dimension two, we treat this case separately in
\S\ref{subsec:lower2}. The proof for this case can basically be
found in~\cite{ego-kon}, and we bring it here for the sake of
clarity and completeness.

For the dimension two case we also
bring a new proof in \S\ref{sec:nazpolsod}. Moreover, this proof shows that
a big inscribed ball can be taken to be with center at a maximal point of the
eigenfunction in the nodal domain.
This proof is due to F.~Nazarov, L.~Polterovich and M.~Sodin and is based on complex
analytic methods.

\subsection{A Short Background}
Related to the problem discussed in this paper is the problem of
estimating the $(n-1)$-Hausdorff measure $H_{n-1}(\lambda)$ of the
nodal set, i.e.~the set where an eigenfunction vanishes.
J.~Br\"uning and D.~Gromes proved in~\cite{brun-gromes}
and~\cite{bruning} sharp lower estimates in dimension two. Namely,
they showed $H_1(\lambda)\geq C\sqrt{\lambda}$. An estimate of the
constant $C$ is given in~\cite{savo}. Later, S.~T.~Yau conjectured
that in any dimension $C_1 \sqrt{\lambda}\geq H_{n-1}(\lambda)\geq
C_2\sqrt{\lambda}$. This was proved in the case of analytic
metrics by H.~Donnelly and C.~Fefferman in~\cite{don-fef-yau}.

 Regarding the inner radius of nodal domains,
 we would like to mention the recent work of B.~Xu \cite{xubin}, in
which he obtains a sharp lower bound on the inner radius for at
least two nodal domains, and the work of V.~Maz'ya and M.~Shubin
\cite{maz-shub}, in which they give sharp bounds on the inner
capacity radius of a nodal domain.

\subsection{Acknowledgements}
I am grateful to Leonid Polterovich for introducing me the problem
and for fruitful discussions. I would like to thank Joseph
Bernstein, Lavi Karp and Mikhail Sodin for enlightening
discussions. I am thankful to Sven Gnutzmann for showing me the
nice nodal domains pictures he generated with his computer program
and for nice discussions. I owe my gratitude also to Moshe Marcus,
Yehuda Pinchover and Itai Shafrir for explaining to me the
subtleties of Sobolev spaces.

I would like to thank Leonid Polterovich, Mikhail Sodin and
F\"edor Nazarov for explaining their proof in dimension two to me,
and for letting me publish it in \S\ref{sec:nazpolsod}.

Special thanks are sent to Sagun Chanillo, Daniel Grieser, Mikhail
Shubin, Bin Xu and the anonymous referee for their comments and corrections
on the first manuscript of this paper.

\section{The Lower Bound on the Inner Radius}
\label{sec:lower}
In this section we prove the existence of the lower bounds on the
inner radius given in Theorems~\ref{thm:inrad}
and~\ref{thm:inrad2}.
\subsection{Lower Bound in Dimension $\geq$ 3}
\label{sec:lower3}
In this section we prove the existence of the lower bound in
Theorem~\ref{thm:inrad}. The proof also gives a bound in the case
where $\dim M = 2$, namely $r_\lambda\geq
C/\sqrt{\lambda\log\lambda}$, but in the next section we treat
this case separately to get rid of the $\sqrt{\log\lambda}$
factor.

 Let $\{\sigma_i\}$ be a finite cellulation of $M$ by cubes, such that for each
$i$ we can put a Euclidean metric $e_i$ on $\sigma_i$, which
satisfies $e_i/4 \leq g \leq 4e_i$.
Let $r_{\lambda,i}$ be the
inner radius of
 $U_{\lambda, i}=U_\lambda\cap \sigma_i$, and $r_{\lambda, i, e}$ be the Euclidean
 inner radius of $U_{\lambda, i}$.  Notice that
 \begin{equation}
   r_{\lambda, i, e}\leq 2r_{\lambda,i} \leq 2r_\lambda\, .
 \end{equation}

\noindent\underline{\scshape Step 1}.
\begin{figure}
$$\quad\quad\includegraphics{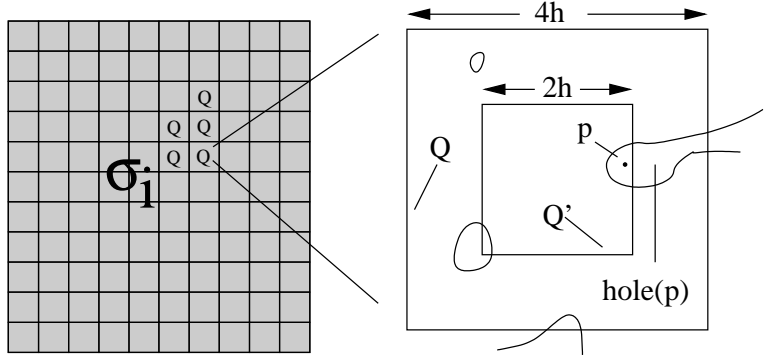}$$
\caption{Proof of Lower Bound on Inner Radius} \label{fig:inrad}
\end{figure}
(See Fig.~\ref{fig:inrad}).
We consider $\sigma_i$ as a compact
 cube in $\RR^n$, with edges parallel to the axes directions.
 We cover $\sigma_i$ by non-overlapping small cubes with edges of size
 $4h$, where $r_{\lambda, i, e}<h<2r_{\lambda, i, e}$.
 Let $Q$ be a copy of one of these small cubes.
 Let $Q'$ be a concentric cube with parallel edges of size~$2h$.

\noindent\underline{\scshape Step 2}. We note that each copy of
$Q'$ contains a point $p\in\sigma_i\setminus U_\lambda$.

Otherwise, we would have $r_{\lambda, i, e}\geq h$, which would
contradict the definition of $h$.

\noindent\underline{\scshape Step 3}. Denote by $\hole(p)$ the
connected component of $Q\setminus (\sigma_i\cap U_\lambda)$ which
contains~$p$. We claim
\begin{equation}
\label{eqn:loc-cour}
 \frac{\Vol_e(\hole(p))}{\Vol_e(Q)} \geq \frac{C_1}{\lambda^{\alpha(n)}(\log\lambda)^{4n}},
\end{equation}
where $\alpha(n)=2n^2+n/4$, and $\Vol_e$ denotes the Euclidean
volume.
 We will denote the right hand side term of~(\ref{eqn:loc-cour}) by~$\gamma(\lambda)$.

Indeed, $\hole(p)$ is a connected component of $U_\lambda'\cap Q$
for some $\lambda$-nodal domain $U_\lambda'$. Hence, we can apply
the following Local Courant's Nodal Domain Theorem.
 \begin{theorem}[\cite{don-fef, chan-mucken, lu}]
  \label{thm:loc-cour}
  Let $B\subseteq M$ be a fixed ball. Let $B'$ be a concentric
  ball of half the radius of $B$. Let $U_\lambda$ be a
  $\lambda$-nodal domain which intersects $B'$.
  Let $B_\lambda$ be a connected component of $B\cap U_\lambda$.
  Then
  \begin{equation}
  \label{eqn:df}
    \Vol(B_\lambda)/\Vol(B)\geq\frac{C_2}{\lambda^{\alpha(n)}(\log\lambda)^{4n}},
  \end{equation}
   where $\alpha(n)=2n^2+n/4$.
\end{theorem}

We remark that in our case~(\ref{eqn:df}) is true also for the
quotient of Euclidean volumes, since the Euclidean metric on
$\sigma_i$ is comparable with the metric coming from $M$.

\noindent\underline{\scshape Step 4}. We let $\tldvphi_{\lambda}
=\chi(U_\lambda)\vphi_{\lambda}$, where $\chi(U_\lambda)$ is the
characteristic function of $U_\lambda$, and similarly,
$\tldvphi_{\lambda, i} = \chi(U_\lambda\cap \sigma_i)
\vphi_\lambda$. Then we have the inequality
\begin{equation}
\label{ineq:loc_lambda} \int_{Q} |\tldvphi_{\lambda,
i}|^2\,\dif(\vol) \leq \beta(\lambda) h^2 \int_Q
|\nabla\tldvphi_{\lambda,i}|^2\,\dif(\vol)\, ,
\end{equation}
where $$\beta(\lambda) = \left\{\begin{array}{lcr}
 C_3\log(1/\gamma(\lambda)) &,& n=2\, , \\
 C_4/(\gamma(\lambda))^{(n-2)/n} &,& n\geq 3\, .
\end{array}\right.$$

\begin{proof}
 Observe that $\tldvphi_{\lambda, i}$ vanishes on
$\hole(p)$. We will use the following Poincar\'{e} type inequality
due to Maz'ya. We discuss it in \S\ref{subsec:poinca-cap}. A
general version of this inequality with weights instead of
Lebesgue measure is proved in~\cite{chan-wheed}.
\begin{theorem}
\label{thm:poinca-vol} Let $Q\subset \RR^n$ be a cube whose edge
is of length~$a$. Let $0<\gamma<1$. Then,
  $$\int_Q |u|^2\,\dif(\vol)\leq
  \beta a^2\int_Q |\nabla u|^2\,\dif(\vol)$$
  for all Lipschitz functions $u$ on $Q$, which vanish on a set
  of measure $\geq\gamma a^n$, and where
  $$\beta = \left\{\begin{array}{lcr}
 C_{5}\log(1/\gamma) &,& n=2\, ,\\
 C_{6}/\gamma^{(n-2)/n} &,& n\geq 3\, .
\end{array}\right.$$
\end{theorem}
 From~(\ref{eqn:loc-cour}) and Theorem~\ref{thm:poinca-vol} applied to $\tldvphi_{\lambda,
 i}$, it follows
\begin{equation*}
\int_{Q} |\tldvphi_{\lambda,i}|^2\,\dif(\vol_e) \leq
\beta(\lambda) h^2 \int_Q
|\nabla_e\tldvphi_{\lambda,i}|_e^2\,\dif(\vol_e) \, .
\end{equation*}
Since the metric on $\sigma_i$ is comparable to the Euclidean
metric, we have also inequality~(\ref{ineq:loc_lambda}).
\end{proof}

\noindent\underline{\scshape Step 5}.
\begin{equation}
\label{ineq:triangle} \int_{\sigma_i}
|\tldvphi_\lambda|^2\,\dif(\vol) \leq 16\beta(\lambda) r_\lambda^2
\int_{\sigma_i} |\nabla\tldvphi_\lambda|^2\,\dif(\vol) \, .
\end{equation}

This is obtained by summing up
inequalities~(\ref{ineq:loc_lambda}) over all cubes~$Q$ which
cover $\sigma_i$, and recalling that $h< 2r_{\lambda, i, e}\leq
4r_\lambda$.

\noindent\underline{\scshape Step 6}. We sum
up~(\ref{ineq:triangle}) over all cubical cells $\sigma_i$ to
obtain a global inequality.
\begin{align}
\label{ineq:global} \int_{U_\lambda}& |\vphi_\lambda|^2
\,\dif(\vol) =
 \int_M |\tldvphi_\lambda|^2 \,\dif(\vol) =
 \sum_i\int_{\sigma_i} |\tldvphi_\lambda|^2\,\dif(\vol) \nonumber\\ & \leq
 16\beta(\lambda) r_\lambda^2 \sum_i \int_{\sigma_i} |\nabla\tldvphi_\lambda|^2 \,\dif(\vol)
 = 16\beta(\lambda) r_\lambda^2 \int_{M}
 |\nabla\tldvphi_\lambda|^2\,\dif(\vol)\nonumber\\
  &= 16\beta(\lambda) r_\lambda^2 \int_{U_\lambda}
 |\nabla\vphi_\lambda|^2\,\dif(\vol)
\end{align}

\noindent\underline{\scshape Step 7}.
$$r_\lambda\geq\left\{
\begin{array}{lcr}
C_{7}/\sqrt{\lambda\log\lambda} &, & n=2\, , \\
C_{8}/\lambda^{n^2-15n/8+1/4}(\log\lambda)^{2n-4}&, & n\geq 3\, .
\end{array}\right.$$

Indeed, by~(\ref{ineq:global})
$$\lambda = \frac{\int_{U_\lambda}
 |\nabla\vphi_\lambda|^2\,\dif(\vol)}
 {\int_{U_\lambda} |\vphi_\lambda|^2  \,\dif(\vol)} \geq
 \frac{1}{16\beta(\lambda)r_\lambda^2}\, .$$

 Thus,
 \begin{align*}
  r_\lambda &\geq \frac{1}{4\sqrt{\lambda\beta(\lambda)}}=
 \left\{\begin{array}{lcr}
   C_{9}/\sqrt{\lambda\log(1/\gamma(\lambda))} &,& n=2\, , \\
   C_{10}\gamma(\lambda)^{(n-2)/2n}/\sqrt{\lambda} &,& n\geq 3\, .
 \end{array}\right. \\
 & \geq
 \left\{\begin{array}{lcr}
  C_{11}/\sqrt{\lambda\log\lambda} &,& n=2\, , \\
  C_{12}/\left(\lambda^{n^2-15n/8+1/4}(\log\lambda)^{2n-4}\right) &,& n\geq 3\, .
 \end{array}\right.
\end{align*}
\qed
\subsection{Lower Bound in Dimension $=$ 2}
\label{subsec:lower2}
We prove the existence of the lower bound on the inner radius in
Theorem~\ref{thm:inrad2}. The arguments below can basically be
found in~\cite{ego-kon}, Chapter~7.

 We begin the proof of Theorem~\ref{thm:inrad2} with Step~1 and
 Step~2 of \S\ref{sec:lower3}. We proceed as follows.

\noindent\underline{\scshape Step 3'}. If $\hole(p)$ does not
touch $\partial Q$
 $$\frac{\Area_e(\hole(p))}{\Area_e(Q)}\geq C_{1}\, ,$$
 where $\Area_e$ denotes the Euclidean area.

\begin{proof}
We recall the Faber-Krahn
inequality in $\RR^n$.
 \begin{theorem}
 \label{thm:faber-krahn}
 Let $\Omega\subseteq\RR^n$ be a bounded domain.
 Then $\lambda_1(\Omega) \geq C_{2}/\Vol(\Omega)^{2/n}$
 \end{theorem}
  We apply Theorem~\ref{thm:faber-krahn} with $\Omega=\hole(p)$.
  We emphasize that $\lambda_1(\hole(p), g) \geq
  C_{3}\lambda_1(\hole(p), e)$,
   since the two metrics are comparable.

  Thus, we obtain
   $$\lambda = \lambda_1(\hole(p), g) \geq \frac{C_{4}}{\Area_e(\hole(p))}\, ,$$
   or, written differently,
$\Area_e(\hole(p))\geq C_{4}/\lambda$. On the other hand,
$\Area_e(Q)=(4h)^2\leq 64 r_\lambda^2\leq 64C_5/\lambda$, where
the last inequality is the upper bound on the inner radius in
Theorem~\ref{thm:inrad2}. So take $C_{1}=C_{4}/(64C_5)$.
\end{proof}

\noindent\underline{{\scshape Step 4'} (part a)}.
  There exists an edge of $Q$, on which the orthogonal projection of
  $\hole(p)$ is of Euclidean size $\geq\gamma\cdot 4h$, where $0<\gamma<1$ is independent
  of $\lambda$.

  Let us denote by $|\pr(\hole(p))|$ the maximal size of the
  projections of $\hole(p)$ on one of the edges of $Q$.
   If $\hole(p)$ touches $\partial Q$, then $|\pr(\hole(p))|\geq 4h/4=h$, and we can take $\gamma=1/4$.
  Otherwise, by Step~3'
  \begin{align*}
  \label{ineq:pr}
   |\pr&(\hole(p))| \geq \sqrt{\Area_e(\hole(p))}\\
   &\geq \sqrt{C_{1} (4h)^2} = 4\sqrt{C_{1}}h\, .
  \end{align*}
  So, we can take $\gamma=\sqrt{C_{1}}$.

 \noindent\underline{{\scshape Step 4'} (part b)}.
\begin{equation}
\label{eqn:locpoinca2}
  \int_Q |\tldvphi_{\lambda,i}|^2\,\dif \vol_e \leq C_{6} h^2 \int_Q
  |\nabla\tldvphi_{\lambda,i}|^2 \,\dif \vol_e .
\end{equation}

  Notice that $\tldvphi_{\lambda, i}$ vanishes on $\hole(p)$.
  Hence, Step~4' (part a) permits us to apply the following Poincar\'{e} type
  inequality to $\tldvphi_{\lambda,i}$. Its proof is given in
  \S\ref{sec:poinca-2}. An inequality in the same spirit can be
  found in~\cite{leonsimon}.
  \begin{theorem}[\cite{ego-kon}, ch.~7]
    \label{thm:poinca-dim2}
     Let $Q\subseteq \RR^2$ be a cube whose edge is of length $a$.
     Let $u$ be a Lipschitz function on $Q$ which vanishes on a curve whose
      projection on one of the edges is of size $\geq \gamma a$. Then
     $$\int_Q |u|^2\,\dif x \leq C(\gamma) a^2 \int_Q |\nabla u|^2\,\dif x\, .$$
\end{theorem}

\noindent\underline{\scshape Steps 5'--7'}.
 To conclude we continue in the same way as in Steps~5--7 of \S\ref{sec:lower3}. \qed
\section{A New Proof in Dimension Two}
\label{sec:nazpolsod}
This section is due to L.~Polterovich, M.~Sodin and F.~Nazarov. In
dimension two we give a proof based on the harmonic measure and
the fact due to Nadirashvili that an eigenfunction on the scale
comparable to the wavelength is almost harmonic in a sense to be defined below.
This proof also gives information about the location of a big ball inscribed in the nodal domain $U_\lambda$.
Namely, we show that if $\phi_\lambda(x_0) = \max_{U_{\lambda}}|\phi_\lambda|$, then
one can find a ball of radius $C/\sqrt{\lambda}$ centered at $x_0$ and inscribed in $U_\lambda$.

Let $D_p\subseteq\Sigma$ be a metric disk centered at $p$.
 Let $f$  be a function defined on $D$. Let $\DD$
 denote the unit disk in $\CC$.
\begin{definition}
We say that $f$ is \emph{$(K, \delta)$-quasi\-harmonic} if there
exists a $K$-quasi\-conformal homeomorphism $h:D\to\DD$, a
harmonic function $u$ on $\DD$, and a function $v$ on $\DD$ with
$1-\delta\leq v\leq 1$, such that
\begin{equation}
  f= (v\cdot u)\circ h\, .
\end{equation}
\end{definition}
\remark. We will assume without loss of generality that $h(p) =0$.

\begin{theorem}[\cite{nadir, naz-pol-sod}]
\label{thm:harmon}
 There exist $K, \eps, \delta>0$ such that for every eigenvalue $\lambda$
 and disk $D\subseteq\Sigma$ of radius $\leq\eps/\sqrt{\lambda}$,
 $\vphi_\lambda|_{D}$ is $(K, \delta)$-quasi\-harmonic.
\end{theorem}

We now choose a preferred system of conformal coordinates on
$(\Sigma, g)$.
\begin{lemma}
\label{lem:coor} There exist positive constants $q_+, q_-, \rho$
such that for each point $p\in M$, there exists a disk
$D_{p,\rho}$ centered at $p$ of radius $\rho$, a conformal map
$\Psi_p:\DD\to D_{p, \rho}$ with $\Psi_p(0)=p$, and a positive
function $q(z)$ on $\DD$ such that
$$\Psi_p^{*}(g) = q(z)|dz|^2, $$
with $q_-<q < q_+$.
\end{lemma}

 Let us take a point $p$, where $|\vphi_\lambda|$ admits its
 maximum on $U_\lambda$. Let $R=\eps/\sqrt{\lambda q_+}$.
  Let $D_{p, R\sqrt{q_+}}\subseteq D_{p, \rho}$ be a disk of radius
  $R\sqrt{q_+}$ centered at $p$.

 We now take the functions $u, v$ defined on $\DD$ which correspond
 to $\vphi_\lambda|_{D_{p, R\sqrt{q_+}}}$ in Theorem~\ref{thm:harmon}.
 We observe that
 \begin{equation*}
  \vphi_\lambda(p)=
  u(0)v(0)\geq u(z)v(z)\geq u(z)(1-\delta) \, ,
 \end{equation*}
 for all $z\in\DD$. Hence
 \begin{equation}
  \label{ineq:max}
   u(0) \geq (1-\delta) \max_\DD u \, .
 \end{equation}
 Now we apply the harmonic measure technic.
 Let $U_{\lambda}^0\subseteq \DD$ be the connected component of
 $\{u>0\}$, which contains $0$. Let $E=\DD\setminus U_{\lambda}^0$.
 Let $\omega$ be the harmonic measure of $E$ in $\DD$. $\omega$ is a bounded harmonic
 function on $U_\lambda^0$, which tends to $1$ on $\partial
 U_\lambda^0\cap \mathrm{Int}(\DD)$ and to $0$ on the interior points
of $\partial U_\lambda^0\cap \partial \DD$.
 Let $r_0=\inf\{|z| :\, z\in E\}$.

 By the Beurling-Nevanlinna theorem (\cite{conf-inv}, sec.~3-3),
 \begin{equation}
 \label{ineq:be-ne-nad}
  \omega(0) \geq 1-C_1\sqrt{r_0}\, .
 \end{equation}
 By the majorization principle
 \begin{equation}
 \label{ineq:major-nad}
  u(0)/\max u \leq 1-\omega(0)\, .
 \end{equation}

 Combining inequalities~(\ref{ineq:max}),~(\ref{ineq:be-ne-nad})
 and~(\ref{ineq:major-nad}) gives us
 \begin{equation}
 \label{ineq:r-c}
  r_0\geq C_2\, .
 \end{equation}

In the final step we apply a distortion theorem proved by Mori for
quasi\-conformal maps. Denote by $\DD_r\subseteq \CC$ the disk
$\{|z| < r\}$.
 Observe that $$\Psi_p(\DD_{R})\subseteq D_{p,R\sqrt{q_+}}\, .$$
 Hence, we can compose
 $$\tldh =h\circ \Psi_p :\DD_{R}\to\DD\, .$$
 $\tldh$ is a $K$-quasi\-conformal map.
 By Mori's Theorem (\cite{quasiconf}, Ch.~III.C) it
 is $\frac{1}{K}$-H\"{o}lder. Moreover, it satisfies an inequality
 \begin{equation}
 \label{ineq:mori}
   |\tldh(z_1)-\tldh(z_2)|\leq M\left(\frac{|z_1-z_2|}{R}\right)^{1/K},
 \end{equation}
  with $M$ depending only on $K$.
 Inequalities~(\ref{ineq:r-c}) and~(\ref{ineq:mori}) imply that
 \begin{equation}
   \frac{\dist(p, \partial (U_\lambda\cap D_{p, R}))}{R}\geq
   \left(\frac{C_2}{M}\right)^K\sqrt{q_-}\, .
 \end{equation}

 Hence,
  \begin{equation}
   \inrad(U_\lambda)\geq
   \left(\frac{C_2}{M}\right)^K\sqrt{q_-} R = C_3/\sqrt{\lambda}\, ,
 \end{equation}
 as desired.
\section{A Review of Poincar\'{e} Type Inequalities}
\label{sec:poinca}
We give an overview of several Poincar\'e type inequalities. In
particular, we prove Theorem~\ref{thm:poinca-vol} and
Theorem~\ref{thm:poinca-dim2}.

\subsection{Poincar\'e Inequality and Capacity}
\label{subsec:poinca-cap}
 Theorem~\ref{thm:poinca-vol} is a direct corollary of
the following two inequalities proved by Maz'ya.
\begin{theorem}[\cite{mazya63}, \S 10.1.2 in \cite{mazya-book}]
\label{thm:poinca-cap} Let $Q\subseteq\RR^n$ be a cube whose edge
is of length~$a$. Let $F\subseteq Q$. Then
  $$\int_{Q} |u|^2\,\dif(\vol)\leq
  \frac{C_1 a^n}{\CAP(F, 2Q)}\int_Q |\nabla u|^2\,\dif(\vol)$$
  for all Lipschitz functions $u$ on $Q$ which vanish on $F$.
\end{theorem}
A few remarks:
\begin{itemize}
\item[(a)]
 $2Q$ denotes a cube concentric with
$Q$, with parallel edges of size twice as large.
\item[(b)] If $\Omega\subseteq\RR^n$ is an open set, and
$\bar{F}\subseteq \Omega$, then $\CAP(F, \Omega)$ denotes the
$L^2$-capacity of $F$ in $\Omega$, namely
 $$\CAP(F, \Omega) = \inf_{u\in\Fl}\left\{\int_\Omega |\nabla u|^2\,\dif x\right\},
 $$
 where $\Fl=\{ u\in C^{\infty}(\Omega),\
 u\equiv 1 \mbox{ on } F,\ \supp (u)\subseteq \Omega\}. $
\item[(c)]
 By Rademacher's Theorem~(\cite{ziemer}), a Lipschitz function is
differentiable almost everywhere, and thus the right hand side has
a meaning.
\item[(d)] A generalization of the inequality to a body which is starlike with respect to a ball
is proved in~\cite{maz-shub-disc}.
\end{itemize}

The next theorem is a capacity--volume inequality.
\begin{theorem}[\S2.2.3 in \cite{mazya-book}]
$$ \CAP (F, \Omega) \geq \left\{\begin{array}{lcr}
C_2/\log (\Area(\Omega)/\Area(F)) &,& n=2\, , \\
C_3 /(\Vol(F)^{-(n-2)/n}-\Vol(\Omega)^{-(n-2)/n})   &,& n\geq 3\,
.
\end{array}\right. $$
In particular, for $n\geq3$ we have
\begin{equation*}
 \CAP(F, \Omega) \geq C_3 \Vol(F)^{(n-2)/n} \, .
\end{equation*}
\end{theorem}

\subsection{A Poincar\'e Inequality in Dimension Two}
\label{sec:poinca-2}
In this section we prove Theorem~\ref{thm:poinca-dim2}. The proof
can be found in chapter~7 of~\cite{ego-kon}. We bring it here for
the sake of clarity.

\begin{proof}
 Let the coordinates be such that $Q=\{0\leq x_1, x_2 \leq a\}$.
 Let the given edge be $Q\cap \{x_1=0\}$, and let $\pr$ denote
 the projection from $Q$ onto this edge.
  Set $E = \pr^{-1}(\pr(\hole(p)))$.
  We claim
    \begin{equation}
     \label{ineq:E-Q}
     \int_{E} |u|^2\,\dif x
      \leq a^{2}
      \int_Q |\nabla u|^2 \,\dif x.
    \end{equation}

  Indeed, let $E_t := E\cap\{x_2 = t\}$.

  We recall the following Poincar\'e type inequality in dimension
  one whose proof is given below.
\begin{lemma}
  \label{ineq:poin-1}
  \begin{equation}
     \label{eqn:poin-1}
     \int_a^b |u|^2\,\dif x\leq
      |b-a|^2\int_a^b |u'|^2\,\dif x
  \end{equation}
  for all Lipschitz functions $u$ on $[a, b]$ which vanish at a point
  of $[a,b]$.
\end{lemma}
  By this lemma
  $$\int_{E_{t}} |u|^2\, \dif x_1
  \leq a^2
  \int_{E_t} |\partial_1 u(x_1, t)|^2 \,\dif x_1.$$
  Integrating over $t\in\pr(\hole(p))$ gives us~(\ref{ineq:E-Q}).

Next we show
\begin{equation}
  \label{ineq:var2}
  \int_Q |u|^2 \, \dif x \leq C_1 a^2
  \int_Q |\nabla u|^2 \,\dif x.
\end{equation}

  By the mean value theorem $\exists t_0$ such that
\begin{equation}
  \label{ineq-avg}
  \int_{E_{t_0}} |u|^2\, \dif x_1 \leq \frac{1}{\gamma\cdot a}
  \int_E |u|^2\,\dif x.
\end{equation}

In addition, we have
\begin{eqnarray*}
|u(x)|^2 &\leq& 2|u(x_1, t_0)|^2
+2 |u(x)-u(x_1, t_0)|^2 \\
&\leq& 2|u(x_1, t_0)|^2
+2 \left(\int_{t_0}^{x_2} |\partial_2 u(x_1, s)|\,\dif s\right)^2 \\
&\leq& 2|u(x_1, t_0)|^2 +2\cdot a\int_{0}^{a} |\partial_2 u(x_1,
s)|^2\,\dif s.
\end{eqnarray*}

Integrating the last inequality over $Q$ gives us
\begin{equation}
\label{ineq-mid}
 \int_Q |u|^2\,\dif x \leq 2\cdot a\int_{E_{t_0}} |u(x_1,
 t_0)|^2\,\dif x_1 + 2a^2\int_Q |\partial_2 u|^2 \,\dif x\, .
\end{equation}

Finally, we combine~(\ref{ineq:E-Q}),~(\ref{ineq-avg})
and~(\ref{ineq-mid}) to get~(\ref{ineq:var2}).
\begin{align*}
  \int_Q & |u|^2\,\dif
x \leq 2 \cdot a\frac{1}{\gamma a}\int_E |u|^2\,\dif x + 2\cdot
a^2 \int_Q |\nabla u|^2\,\dif x\\ & \leq C_1 a^2
 \int_Q |\nabla u|^2 \,\dif x\,.
\end{align*}
\end{proof}
\subsection{A Poincar\'{e} Inequality in Dimension One}
We prove Lemma~\ref{ineq:poin-1}.
\begin{proof}
  By scaling, it is enough to prove~(\ref{eqn:poin-1}) for the segment $[0,1]$.
  Suppose $u(x_0)=0$. Since a Lipschitz function is
  absolutely continuous, we have

  $$|u(x)|^2
  =\left|\int_{x_0}^x u'(t)\,\dif t\right|^2\leq \int_0^1 \left|u'(t)\right|^2\,\dif t.$$

  We integrate over $[0,1]$ to get the desired inequality.
\end{proof}

\section{$\lambda_1$ and Inner Radius}
\label{sec:gen-bounds}
We prove Theorem~\ref{thm:inrad-gen}, which relates the inner
radius to $\lambda_1$.
\begin{proof}
  Let $\{V_i\}$ be a finite open cover of $M$, such that
  for each $i$ one can put a Euclidean metric $e_i$ on $V_i$, which
  satisfies $e_i/4\leq g\leq 4e_i$.
  Let $\alpha$ be the Lebesgue number of the covering.

  Let $r=\min(\inrad(\Omega), \alpha)$.
  Let $B\subseteq \Omega$ be a ball of radius $r$.
  We can assume that $B\subseteq V_1$.
  Let $B_e\subseteq B$ be a Euclidean ball of radius $r/2$.
  By monotonicity of $\lambda_1$, we know that $\lambda_1(\Omega, g)\leq
  \lambda_1(B, g)\leq
  \lambda_1(B_e, g)$, but since the Riemannian metric
  on $B_e$ is comparable to the Euclidean metric on it, it follows
  from the variational principle that
   $$\lambda_1(B_e, g) \leq C_1\lambda_1(B_e, e_1) = C_2/r^2\leq
   C_3/\inrad(\Omega)^2,$$
where in the last inequality we used the fact that $\inrad(\Omega)\leq C_4\alpha$.
\end{proof}

\remark. We would like to emphasize that in general there is no
lower bound on $\lambda_1$ is terms of the inner radius. However,
in dimension two, as pointed out to us by Daniel Grieser and
Mikhail Shubin, there exists a lower bound on $\lambda_1$ in terms
of the inner radius and the connectivity of $\Omega$. This was
proved in~\cite{hayman} and~\cite{osserman}. For a more detailed
account of the subject one can consult~\cite{lieblow}.

\providecommand{\bysame}{\leavevmode\hbox to3em{\hrulefill}\thinspace}
\providecommand{\MR}{\relax\ifhmode\unskip\space\fi MR }
\providecommand{\MRhref}[2]{%
  \href{http://www.ams.org/mathscinet-getitem?mr=#1}{#2}
}
\providecommand{\href}[2]{#2}

\noindent Dan Mangoubi,\\
Department of Mathematics,\\
The Technion, \\
Haifa 32000,\\
ISRAEL.

\noindent email: mangoubi@techunix.technion.ac.il
\end{document}